\documentclass{amsart}
\usepackage{amsthm}                 
\usepackage{amsmath}                
\usepackage{amsfonts,amssymb}       
\usepackage{graphicx}               
\usepackage{float}                  
\usepackage{dsfont}
\usepackage{rotating}               
\restylefloat{figure}               
\usepackage{mathrsfs}               
\usepackage{fancyhdr}               
\usepackage[breaklinks]{hyperref}   
\usepackage[numbers,sort&compress]{natbib}
\theoremstyle{plain}
\newtheorem{lem}{Lemma}
\newtheorem{prop}[lem]{Proposition}
\newtheorem{thm}[lem]{Theorem}
\newtheorem{cor}[lem]{Corollary}

\theoremstyle{definition}

\newtheorem{exmp}[lem]{Example}

\newcommand{\R}{\mathbb R}

\newcommand{\Z}{\mathbb Z}
\newcommand{\N}{\mathbb N}

\newcommand{\Diff}{\mbox{\rm Diff}}

\renewcommand{\L}{\mathcal L}
\newcommand{\A}{\mathbb A}
\newcommand{\id}{\text{\rm id}}
\newcommand{\dx}{\,\text{\rm d}x}
\renewcommand{\d}{\,\text{\rm d}}
\newcommand{\dw}{\text{\rm d}}

\newcommand{\ad}{\text{\rm ad}}

\renewcommand{\S}{\mathbb S}

\renewcommand{\phi}{\varphi}

\newcommand{\norm}[1]{\left|\!\left|#1\right|\!\right|}

\newcommand{\ska}[2]{\left\langle #1,#2\right\rangle}

\renewcommand{\i}{\text{i}}

\newcommand{\bea}{\begin{eqnarray}}
\newcommand{\eea}{\end{eqnarray}}
\newcommand{\beq}{\begin{equation}}
\newcommand{\eeq}{\end{equation}}
\renewcommand{\phi}{\varphi}

\renewcommand{\autoref}[1]{\text{Eq.}~\eqref{#1}}
\begin{document}
\title{A note on multi-dimensional Camassa-Holm type systems on the torus}
\author{Martin Kohlmann}
\address{Peter L. Reichertz Institute for Medical Informatics, University of Braunschweig, D-38106 Braunschweig, Germany}
\email{martin.kohlmann@plri.de}
\begin{abstract} We present a $2n$-component nonlinear evolutionary PDE which includes the $n$-dimensional versions of the Camassa-Holm and the Hunter-Saxton systems as well as their partially averaged variations. Our goal is to apply Arnold's \cite{A66,EM70} geometric formalism to this general equation in order to obtain results on well-posedness, conservation laws or stability of its solutions. Following the line of arguments of the paper \cite{K11} we present geometric aspects of a two-dimensional periodic $\mu$-$b$-equation on the diffeomorphism group of the torus in this context.
\end{abstract}
\maketitle
\tableofcontents
\section{Introduction}
In this paper, we study the following system of nonlinear partial differential equations
\bea\label{*}
\left\{
  \begin{array}{rcl}
    m_t    & = & -u\cdot\nabla m-(\nabla u)^T\cdot m-m(\nabla\cdot u)-(\nabla\rho)^T\cdot\rho, \\
    \rho_t & = & -\nabla\rho\cdot u-\rho(\nabla\cdot u), \\
  \end{array}
\right.
\eea
where $u$, $m$ and $\rho$ are vector fields on the $n$-torus $\S^n\simeq\R^n/\Z^n$ which also depend on time $t\geq 0$. Furthermore, it is assumed that there is a linear operator $A$ such that $m=Au$, and that $A$ is of the form $\alpha\mu+\beta-\Delta$, with $\alpha,\beta\in\{0,1\}$ and $\alpha+\beta\neq 2$; here $\mu(u)=\int_{\S^n}u(x)\d^nx$ is the mean value operator.

For $(\alpha,\beta)=(0,1)$ and $(\alpha,\beta)=(1,0)$, the system \eqref{*} reduces to the $2n$-component $n$-dimensional periodic $(\mu)-$Camassa-Holm equation. For the choice $\alpha=\beta=0$, one obtains from \eqref{*} the $2n$-component $n$-dimensional periodic Hunter-Saxton system. If in addition $\rho\equiv 0$, the above mentioned systems reduce to the corresponding $n$-component versions. For $n=1$, there are several studies of the $(\mu)-$Camassa-Holm and the Hunter-Saxton equations or systems, respectively.
The Camassa-Holm equation
\beq\label{CH}u_t+3uu_x=2u_xu_{xx}+uu_{xxx}+u_{txx}\eeq
first appeared in 1993 and was introduced by bi-Hamiltonian methods in \cite{CH93}. In the subsequent years, various interesting properties of this equation and its solutions, e.g., integrability \cite{I05,I07,J03}, blow-up and peakons \cite{CE98,M04} and geometric aspects \cite{CK02,CK03,K99,M98,L07}, have been discussed. Related issues have been worked out for the two-component generalization of \eqref{CH}
\bea
\label{2CH}
\left\{
  \begin{array}{rcl}
    m_t    & = & -m_xu-2u_xm-\rho\rho_x, \\
    \rho_t & = & -(\rho u)_x, \\
  \end{array}
\right.
\eea
with $m=u-u_{xx}$, in \cite{CI08,EKL11,GO06,GZ10,HI11,HT09,LL09}. In \cite{KLM08}, the authors proposed the equation
\beq\label{muCH}2\mu(u)u_x=2u_xu_{xx}+uu_{xxx}+u_{txx}\eeq
and called it the $\mu$-Hunter-Saxton equation. This equation also appeared in \cite{LMT10} and has been called the $\mu$-Camassa-Holm equation there. The associated two-component version, which is the system \eqref{2CH} with $m=\mu(u)-u_{xx}$ , has been the subject of \cite{LY11,LY11b,Z10} where the authors establish its bi-Hamiltonian and variational nature as well as blow-up, global existence and the existence of weak solutions. The Hunter-Saxton equation
\beq\label{HS}0=2u_xu_{xx}+uu_{xxx}+u_{txx}\eeq
and its associated two-component analog, which is \autoref{2CH} with $m=-u_{xx}$, also have several interesting mathematical properties, c.f., e.g., \cite{HZ94,L07',L07'',L08,LL09,WW10,W10,W11,Y04}.
Note that the Camassa-Holm and Hunter-Saxton equations have meaningful physical interpretations since they are related to fluid dynamics: In \cite{CL09} it is shown that \autoref{CH} describes the motion of shallow water waves over a flat bottom under the action of gravity. In \cite{HS91} it is explained that \autoref{HS} is derived from the least action principle for orientation waves in liquid crystals.

Counter to the large amount of papers referring to the case $n=1$, higher dimensional variations of the Camassa-Holm type systems named above have rarely been studied. However, multi-variable extensions of these equations are of interest from the physical and the mathematical point of view as explained in, e.g., \cite{HM05,GB09,LMT10}.

Concerning $n=2$, $(\alpha,\beta)=(0,1)$ and $\rho=0$, we have shown in \cite{K11} that the two-dimensional Camassa-Holm equation \cite{KSD01,GB09,HM05} can be interpreted in Arnold's geometric framework, i.e., as a geodesic flow on the diffeomorphism group of the torus. The key idea is to recast the two-dimensional Camassa-Holm equation as a geodesic equation on the torus diffeomorphism group with respect to a suitable right-invariant weak Riemannian metric. In addition, there is a smooth affine connection which preserves the metric. As a consequence of the geometric approach, one obtains results on, e.g., well-posedness and stability.

In fact, many results obtained in \cite{K11} also apply to the general system \eqref{*} which includes the Camassa-Holm and Hunter-Saxton systems in any dimension $n$. The geometric theory is not only aesthetically appealing, but also results in information about solutions or conservation laws. It is the aim of this note to prove that \eqref{*} is well-posed on a scale of Sobolev spaces and also in the smooth category. Concerning the special choice $n=2$, $(\alpha,\beta)=(1,0)$ and $\rho=0$, we consider a two-dimensional periodic $\mu$-$b$-equation \cite{HS03,LMT10} and prove that $b=2$ is the only case for which one obtains a metric Euler equation (namely the 2D $\mu$-Camassa-Holm equation). We also provide some explicit computations of the sectional curvature of the associated torus diffeomorphism group in this case.
\section{The geometric formalism}\label{sec_geo}
We let $\alpha,\beta,\gamma\in\{0,1\}$, $\alpha+\beta\neq 2$, $\mathds{1}_n$ the $n\times n$ identity matrix, and set
$$\A=(1-\gamma)(\alpha\mu+\beta-\Delta)\mathds{1}_n+\gamma\,\text{\rm diag}(\alpha\mu+\beta-\Delta,1)\otimes\mathds{1}_n.$$
For $\gamma=0$, $(\alpha,\beta)\neq(0,0)$ and $s>n/2+1$, let $G^s=\Diff^s(\S^n)$ denote the group of orientation preserving diffeomorphisms $\S^n\to\S^n$ which are elements of the Sobolev space $H^s(\S^n)$. For $\alpha=\beta=\gamma=0$ we let $G^s=\Diff_0^s(\S^n)$ be the subgroup of $\Diff^s(\S^n)$ defined by the condition $p(0)=0$. For $\gamma=1$ and $s>n/2+1$, let $G^s$ denote the semidirect product obtained from the orientation preserving $H^s$ diffeomorphisms of $\S^n$ (or its subgroup $\Diff_0^s(\S^n)$ in the case $\alpha=\beta=0$) with the space $H^{s-1}(\S^n;\R^n)$. Then $G^s$ is a smooth Hilbert manifold and a topological group and
$$T_eG^s\simeq
\left\{
\begin{array}{cc}
  H^s(\S^n;\R^n), & \gamma=0,\\
  H^s(\S^n;\R^n)\times H^{s-1}(\S^n;\R^n), & \gamma=1,
\end{array}
\right.
$$
for $(\alpha,\beta)\neq(0,0)$, and
$$T_eG^s\simeq
\left\{
\begin{array}{cc}
  H^s_0(\S^n;\R^n), & \gamma=0,\\
  H^s_0(\S^n;\R^n)\times H^{s-1}(\S^n;\R^n), & \gamma=1,
\end{array}
\right.
$$
if $\alpha=\beta=0$; here $H^s_0(\S^n;\R^n)$ is the subspace of $H^s(\S^n;\R^n)$ defined by the condition $u(0)=0$.
For $s\to\infty$, we arrive at the smooth category which is described as above with $H^s$ replaced by $C^\infty$. The group $G^\infty=\cap_{s>n/2+1}G^s$ is a Lie group and a smooth Fr\'echet manifold. Further details about the groups under discussion can be found in Sect.~2 of \cite{K11}.

Our strategy is to define a right-invariant metric on $G^s$ such that the corresponding Euler equation is the $n$-dimensional version of the equation specified in the following tabular.
$$
\begin{array}{|l||c|c|c|}
  \hline
  (\alpha,\beta) & (0,0) & (0,1) & (1,0) \\
  \hline
  \gamma=0 & \text{HS}  & \text{CH} & \mu-\text{CH} \\
  \gamma=1 & 2\text{HS} & 2\text{CH} & \mu-2\text{CH}  \\
  \hline
\end{array}
$$
In the framework of Arnold's \cite{A66} approach, we have to proceed as follows:
\begin{enumerate}
\item Check that the \emph{inertia operator} $\A$ is indeed a topological isomorphism $T_eG^\infty\to T_eG^\infty$.
\item Define a right-invariant weak Riemannian metric on $G^\infty$ which coincides with the bilinear form induced by $\A$ at the identity.
\item Define an affine connection $\bar\nabla$ on $G^\infty$ and a corresponding right-invariant \emph{Christoffel operator} $\Gamma_p$ such that the equation $p_{tt}=\Gamma_p(p_t,p_t)$ is the geodesic equation for the Lagrangian flow $p(t)\subset G^\infty$ defined by $\bar\nabla$.
\item Check that the metric and the connection produce indeed the same geodesic flow. The geodesics on $G^\infty$ are then also determined by the solutions of the equation $w_t=\ad^*_ww$, in terms of the Eulerian variable $w=p_t\circ p_1^{-1}$, where $p=p_1$ if $\gamma=0$ or $p=(p_1,p_2)$ if $\gamma=1$.
\end{enumerate}
The geometric approach will result in information about the solutions of the equations under consideration, as we will elucidate in the following. For the sake of simplicity, we will first work with the Sobolev spaces $H^s$ and we will discuss the limit $s\to\infty$ later on.

We begin with the inertia operator $\A$. Let $\hat H^s(\S^n;\R^n)$ denote the subspace of $H^s(\S^n;\R^n)$ defined by the condition $\mu(u)=0$. To keep the notation as simple as possible, we write $H^s(\S^n;\R^n)=H^s$ etc. henceforth. We assume that $s>n/2+2$.
\lem The operator $\A$ is a topological isomorphism
\begin{align}
&H^s\to H^{s-2}, & \gamma=0,(\alpha,\beta)\neq(0,0),\nonumber\\
&H^s_0\to\hat H^{s-2}, & \gamma=0,(\alpha,\beta)=(0,0),\nonumber\\
&H^s\times H^{s-1}\to H^{s-2}\times H^{s-1}, & \gamma=1,(\alpha,\beta)\neq(0,0),\nonumber\\
&H^s_0\times H^{s-1}\to\hat H^{s-2}\times H^{s-1}, & \gamma=1,(\alpha,\beta)=(0,0).\nonumber
\end{align}
\endlem\rm
\proof It suffices to show that $1-\Delta,\mu-\Delta\colon H^s\to H^{s-2}$ and $-\Delta\colon H^s_0\to\hat H^{s-2}$ are topological isomorphisms. For $1-\Delta$ and $\mu-\Delta$ this has been established as explained in \cite{LMT10}. We now prove that the operator $A=-\Delta$ is a topological isomorphism $D(A)=H^s_0\to\hat H^{n-2}$; here, we write $Au=(Au_1,\ldots,Au_n)$ for simplicity. Let $\mathcal F\colon L_2(\S^n)\to\ell_2(\Z^n)$ be the Fourier transform. Any $u\in D(A)$ can be written as $u=\sum_{\alpha\in\Z^n}u_\alpha e^{2\pi\i\alpha\cdot x}$, with $(u_\alpha)_{\alpha\in\Z^n}=\mathcal F(u)$. Since $\ker A=\{0\}$ it follows that $A$ is injective. To see that $A$ is surjective, we pick $v=\sum_{\alpha\in\Z^n\backslash\{0\}}v_\alpha e^{2\pi\i\alpha\cdot x}\in\hat H^{s-2}$ and define
$$w=\sum_{\alpha\in\Z^n\backslash\{0\}}\frac{v_\alpha}{4\pi^2\alpha^2}e^{2\pi\i\alpha\cdot x}-\sum_{\alpha\in\Z^n\backslash\{0\}}\frac{v_\alpha}{4\pi^2\alpha^2}.$$
Since $w(0)=0$ and since $A$ is a Fourier multiplication operator with the symbol $4\pi^2\alpha^2$, it follows that $Aw=v$. Finally, it is easy to derive the estimates $\norm{Au}_{H^{s-2}}\lesssim\norm{u}_{H^{s}}$ and $\norm{A^{-1}v}_{H^{s}}\lesssim\norm{v}_{H^{s-2}}$, for $u\in D(A)$ and $v\in D(A^{-1})$ respectively.
\endproof
The second step is to introduce a suitable right-invariant metric.
\lem The map
\beq\label{skapro}\ska{\cdot}{\cdot}\colon T_eG^s\times T_eG^s\to\R,\quad (u,v)\mapsto\int_{\S^n} u\A v\d^n x\eeq
is a scalar product on $T_eG^s$.
\endlem
\proof Obviously, $\ska{\cdot}{\cdot}$ is a symmetric bilinear form on $T_eG^s$. Since we have
$$
\ska{u}{v}=\sum_{i=1}^n\left[\alpha\mu(u_i)\mu(v_i)+\int_{\S^n}(\beta u_iv_i+\nabla u_i\cdot\nabla v_i)\d^n x\right]
+\gamma\sum_{i=n+1}^{2n}\int_{\S^n}u_iv_i\d^n x\nonumber
$$
it is immediate to check that, for all possible combinations of $\alpha,\beta,\gamma$, we have that $\ska{u}{u}\geq 0$ and $\ska{u}{u}=0$ iff $u=0$.
\endproof
We extend the map defined in \eqref{skapro} to a map
\beq\label{skaprox} G^s\times TG^s\times TG^s\to\R,\quad (p,u,v)\mapsto\ska{u}{v}_p=\ska{u\circ p_1^{-1}}{v\circ p_1^{-1}},\eeq
where $u,v\in T_pG^s$ and $p=p_1$, for $\gamma=0$, and $p=(p_1,p_2)$, for $\gamma=1$. Note that $u\circ p_1^{-1}$ is the differential at the point $p$ of the right translation map $R_{p^{-1}}\colon G^s\to G^s$, $q\mapsto qp^{-1}$.
\lem The pair $(G^s,\ska{\cdot}{\cdot})$ is a weak Riemannian manifold. In particular, the map
$$G^s\to\mathcal L^2_{\text{\rm sym}}(TG^s;\R),\quad p\mapsto\ska{\cdot}{\cdot}_p$$
is smooth.
\endlem\rm
\proof
That $\ska{\cdot}{\cdot}_p$ is a weak Riemannian metric follows as in \P 9 of \cite{EM70}; cf.~Theorem~9.1. That the map $p\mapsto\ska{\cdot}{\cdot}_p$ is smooth also follows from the explicit representation
\begin{align}
\ska{u}{v}_p&=\sum_{i=1}^n\bigg[\alpha\mu(u_i|\nabla p_1|)\mu(v_i|\nabla p_1|)+\beta\int_{\S^n}u_iv_i|\nabla p_1|\d^n x\nonumber\\
&\hspace{-1cm}+\int_{\S^n}[\nabla u_i^T\cdot(\nabla p_1)^{-1}]\cdot[\nabla v_i^T\cdot(\nabla p_1)^{-1}]|\nabla p_1|\d^n x\bigg]
+\gamma\sum_{i=n+1}^{2n}\int_{\S^n}u_iv_i|\nabla p_1|\d^n x\nonumber
\end{align}
and this achieves the proof of the lemma.
\endproof
Now let $A=\alpha\mu+\beta-\Delta$,
\beq
\Gamma_{(\id,0)}((u,\rho),(v,\eta))=
\left(
  \begin{array}{c}
    \Gamma^0_\id(u,v) \\
    0 \\
  \end{array}
\right)-
\frac{1}{2}
\left(
  \begin{array}{c}
   A^{-1}[\nabla(\rho\cdot\eta)] \\
    \rho(\nabla\cdot v)+\eta(\nabla\cdot u) \\
  \end{array}
\right),
\label{Christoffel2}\eeq
where
\begin{align} \Gamma^0_\id(u,v) &= -\frac{1}{2}A^{-1}\{u\cdot\nabla(Av)+(\nabla u)^TA v+Av(\nabla\cdot u)-A(\nabla u\cdot v)\nonumber\\
&\qquad +v\cdot\nabla(Au)+(\nabla v)^T Au+Au(\nabla\cdot v)-A(\nabla v\cdot u)\},
\label{Christoffel1}\end{align}
and define, for $\gamma=0$, $p=p_1\in G^s$, and $u,v\in T_pG^s\simeq H^s$ ($\simeq H^s_0$ respectively),
$$\Gamma_{p}(u,v)=\Gamma_{\id}^0(u\circ p_1^{-1},v\circ p_1^{-1})\circ p_1,$$
and for $\gamma=1$, $p=(p_1,p_2)\in G^s$, and $u,v\in T_pG^s\simeq H^s\times H^{s-1}$ ($\simeq H^s_0\times H^{s-1}$ respectively),
$$\Gamma_{p}(u,v)=\Gamma_{(\id,0)}(u\circ p_1^{-1},v\circ p_1^{-1})\circ p_1.$$
It is easy to check that, for $\gamma=0$, $\Gamma_p\in\L^2_{\text{sym}}(H^s,H^s)$ ($\Gamma_p\in\L^2_{\text{sym}}(H^s_0,H^s_0)$, respectively) and, for $\gamma=1$, $\Gamma_p\in\L^2_{\text{sym}}(H^s\times H^{s-1},H^s\times H^{s-1})$ ($\Gamma_p\in\L^2_{\text{sym}}(H^s_0\times H^{s-1},H^s_0\times H^{s-1})$, respectively).

For a given function $w\in T_eG^s$, there is a unique solution to the initial value problem
$$
\left\{
\begin{array}{rcll}
  p_t(t) & = & w\circ p_1(t), \\
  p(0) & = & e
\end{array}
\right.
$$
on the Hilbert manifold $G^s$. This gives us the \emph{Lagrangian coordinates} $p=p(t)$ for which we consider the second order equation
\beq p_{tt}=\Gamma_{p}(p_t,p_t).\label{geoeq}\eeq
\autoref{geoeq} is the geodesic equation on $G^s$ for the affine connection
\beq\label{connection}(\bar\nabla_XY)(p)=DY(p)\cdot X(p)-\Gamma_{p}(X(p),Y(p)),\eeq
where $X,Y$ are smooth vector fields on $G^s$. It remains to prove the following proposition.
\prop Let $\mathfrak X(G^s)$ be the space of smooth vector fields on $G^s$. The map $\bar\nabla\colon\mathfrak X(G^s)\times\mathfrak X(G^s)\to\mathfrak X(G^s)$ defined in \eqref{connection} is a smooth, torsion-free affine connection on $G^s$, i.e.,
\begin{itemize}
\item[(i)] $\bar\nabla_{fX+gY}Z=f\bar\nabla_XZ+g\bar\nabla_YZ$,
\item[(ii)] $\bar\nabla_X(Y+Z)=\bar\nabla_XY+\bar\nabla_XZ$,
\item[(iii)] $\bar\nabla_X(fY)=f\bar\nabla_XY+X(f)Y$,
\item[(iv)] $\bar\nabla_XY-\bar\nabla_YX=[X,Y]$,
\end{itemize}
for all $X,Y,Z\in\mathfrak X(G^s)$ and all $f,g\in C^{\infty}(G^s;\R)$. Moreover, the map $p\mapsto(p,(\bar\nabla_XY)(p))$, $G^s\mapsto TG^s$, is smooth for any $X,Y\in\mathfrak{X}(G^s)$. Finally, the connection $\bar\nabla$ preserves the metric $\ska{\cdot}{\cdot}$ in the usual sense
$$X\ska{Y}{Z}=\ska{\bar\nabla_XY}{Z}+\ska{\bar\nabla_XZ}{Y},\quad\forall X,Y,Z\in\mathfrak{X}(G^s).$$
\endprop\rm
\proof That $\bar\nabla$ satisfies the properties $\text{(i)}-\text{(iv)}$ follows immediately from the definition \eqref{connection} and the local formula $[X,Y]=DY\cdot X-DX\cdot Y$. That $p\mapsto(p,(\bar\nabla_XY)(p))$ is smooth follows from the fact that the geodesic spray $p\mapsto\Gamma_p$ is smooth; this can be deduced as in Appendix~A.1 of \cite{K11}.
If $\gamma=0$, the compatibility condition can be obtained by the line of arguments in the proof of Proposition~5 of \cite{K11}, since the calculations there do not depend on the particular choice of the inertia operator and can be generalized immediately to dimensions greater than two. For $\gamma=1$, it remains to prove that the sum of the following six terms coming from the additional $n$ components
\begin{align}
&-\int_{\S^n}\tilde\Gamma_1(u,v)A w_1\d^nx-\int_{\S^n}\tilde\Gamma_1(u,w)A v_1\d^nx-\int_{\S^n}\Gamma_2(u,v)w_2\d^nx\nonumber\\
&-\int_{\S^n}\Gamma_2(u,w)v_2\d^nx+\int_{\S^n}(\nabla v_2\cdot u_1)w_2\d^nx+\int_{\S^n}(\nabla w_2\cdot u_1)v_2\d^nx\nonumber
\end{align}
is zero; here, $u=(u_1,u_2)$ etc., $\tilde\Gamma_1(u,v)=-\frac{1}{2}A^{-1}[(\nabla u_2)^Tv_2+(\nabla v_2)^Tu_2]$ and $\Gamma_2$ is the second component of the Christoffel map \eqref{Christoffel2}. Using the multi-dimensional integration by parts formula, a careful observation then shows that the above terms indeed cancel out. Thus compatibility of $\bar\nabla$ and $\ska{\cdot}{\cdot}$ holds true for any possible choice of the parameters $\alpha,\beta,\gamma$.
\endproof
%
To recover the Camassa-Holm type systems in terms of the \emph{Eulerian coordinates} $u=u(t)$ or $(u,\rho)=(u,\rho)(t)$ respectively from the geometric picture, we compute the dual operator of the commutator bracket $[\cdot,\cdot]$ with respect to the right-invariant metric induced by $\A$. Using integration by parts one easily verifies that
$$\ska{[u,w]}{v}=\ska{B(u,v)}{w},\quad u,v,w\in T_eG^s,$$
where
$$B(u,v)=-A^{-1}(u\cdot\nabla(A v)+(\nabla u)^TA v+A v(\nabla\cdot u)),$$
for $\gamma=0$, and
$$B(u,v)=
\left(
  \begin{array}{c}
    \!\!-A^{-1}(u_1\cdot\nabla(A v_1)+(\nabla u_1)^TA v_1+A v_1(\nabla\cdot u_1)+(\nabla u_2)^Tv_2)\!\! \\
    -\nabla v_2\cdot u_1-(\nabla\cdot u_1)v_2 \\
  \end{array}
\right),
$$
for $\gamma=1$. It is clear that we can rewrite the systems under discussion in the form $u_t=B(u,u)$ and $(u_t,\rho_t)=B((u,\rho),(u,\rho))$ respectively and that the right-hand side $B$ looks as specified above.
\section{Applications of the geometric theory}
We let $X^s=X^s_{\alpha,\beta,\gamma}$ be the space
$$
\begin{array}{|l||c|c|c|}
  \hline
  (\alpha,\beta) & (0,0) & (0,1) & (1,0) \\
  \hline
  \gamma=0 & H_0^s  & H^s & H^s \\
  \gamma=1 & H_0^s\times H^{s-1} & H^s\times H^{s-1} & H^s\times H^{s-1}  \\
  \hline
\end{array}
$$
and $X^\infty$ the corresponding space in the smooth category.
Since the existence of a smooth connection immediately implies the existence of a smooth geodesic flow, we obtain well-posedness of the geodesic equation \eqref{geoeq} in the $H^s$ category.
\prop Let $\Gamma$ be the Christoffel map defined in Sect.~\ref{sec_geo}, for some admissible triplet $(\alpha,\beta,\gamma)$. There is an open neighborhood $U\subset X^s$ containing zero, for $s>n/2+2$, $n\geq 1$, such that for any $w_0\in U$, the Cauchy problem
\bea\label{IVP_flow}
\left\{
\begin{array}{rcl}
p_{tt} & = & \Gamma_{p}(p_t,p_t), \\
p_t(0)  & = & w_0,\\
p(0)      & = & e
\end{array}
\right.
\eea
for the geodesic flow corresponding to the Camassa-Holm type system specified by $(\alpha,\beta,\gamma)$ on $G^{s}$ has a unique solution $p(t)$, one some interval $[0,T)$, which depends smoothly on time and on the initial value, i.e., the map $(t,w_0)\mapsto p(t)$, $[0,T)\times U\to G^s$ is smooth.
\endprop\rm
Indeed, well-posedness also holds in the smooth category. This follows from an application of Theorem~12.1 in \cite{EM70} where it is shown that the geodesic flow preserves its spatial regularity as we increase the regularity of the initial datum.
\thm\label{lwp_geo} Let $\Gamma$ be the Christoffel map defined in Sect.~\ref{sec_geo}, for some admissible triplet $(\alpha,\beta,\gamma)$. There is an open neighborhood $U\subset X^s$ containing zero, for some fixed $s>n/2+2$, $n\geq 1$, such that for any smooth $w_0\in U$, the Cauchy problem \eqref{IVP_flow} for the geodesic flow corresponding to the Camassa-Holm type system specified by $(\alpha,\beta,\gamma)$ on $G^{\infty}$ has a unique solution $p(t)$, one some interval $[0,T)$, which depends smoothly on time and on the initial value, i.e., the map $(t,w_0)\mapsto p(t)$, $[0,T)\times U_\infty\to G^\infty$ is smooth; here $U_\infty$ denotes the subset of smooth functions lying in $U$.
\endprop\rm
In view of the group properties of $G^s$ and $G^\infty$, we immediately deduce the following well-posedness results for our Camassa-Holm type systems in Euclidean variables.
\cor There exists an open neighborhood $U\subset X^s$ containing zero, $s>n/2+2$, $n\geq 1$, such that for any $w_0\in U$ there exists an interval $J=[0,T)$ and a unique solution
$$w\in C(J,X^s)\times C^1(J,X^{s-1})$$
to the Camassa-Holm type system specified by $(\alpha,\beta,\gamma)$, satisfying the initial condition $w(0)=w_0$, and depending continuously on $w_0$. If $w_0\in U_\infty$, as defined in Theorem~\ref{lwp_geo}, we have that
$$w\in C^{\infty}(J,X^{\infty}),$$
and the local flow $(t,w_0)\mapsto w$ is smooth.
\endcor\rm
Another important consequence of the geometric formalism presented in Sect.~\ref{sec_geo} is the following conservation law. In comparison with the geodesic motion of a three-dimensional rigid body on the finite dimensional Lie group $SO(3)$, we have conservation of the angular momentum in the frame of reference of the rotating body.
\thm\label{lem_conservation} Any solution to the Camassa-Holm type system specified by $\alpha,\beta,\gamma=1$ satisfies the conservation law
\beq\label{consv1}\frac{\dw}{\dw t}
\left(
  \begin{array}{c}
    (\nabla p_1)^T(m\circ p_1)|\nabla p_1|+(\nabla p_2)^T(\rho\circ p_1)|\nabla p_1| \\
    (\rho\circ p_1)|\nabla p_1|\\
  \end{array}
\right)
=0.
\eeq
In addition, for $\alpha=\beta=0$, we find conservation of the quantity
\beq\label{consv2}\frac{\dw}{\dw t}\int_{\S^n}(|\nabla u_1|^2+\cdots+|\nabla u_n|^2+|\rho|^2)\d^n x=0.\eeq
The corresponding conservation laws for $\gamma=0$ are obtained by setting $\rho=0$ in \eqref{consv1} and \eqref{consv2}.
\endthm\rm
\proof Let $\gamma=1$. To prove the conservation law \eqref{consv2}, we multiply the first row equation of the $n$-dimensional Hunter-Saxton system \eqref{*} by $u$, integrate over the $n$-torus and use the multi-dimensional integration by parts formula to obtain that
$$\sum_{i=1}^n\int_{\S^n}\nabla u_{t,i}\cdot\nabla u_i\d^nx=-\int_{\S^n}[(\nabla\rho)^T\cdot\rho]\cdot u\d^nx.$$
Multiplying the second row equation of the Hunter-Saxton system by $\rho$ and integrating over $\S^n$ we find that
$$\int_{\S^n}\rho_t\cdot\rho\d^nx=-\int_{\S^n}(\nabla\rho\cdot u+(\nabla\cdot u)\rho)\cdot\rho\d^nx.$$
Performing integration by parts once again shows that
$$\frac{1}{2}\frac{\dw}{\dw t}\int_{\S^n}(|\nabla u_1|^2+\cdots+|\nabla u_n|^2+|\rho|^2)\d^n x=\int_{\S^n}\left(\sum_{i=1}^n\nabla u_{t,i}\cdot\nabla u_i+\rho_t\cdot\rho\right)\dw^n x=0.$$
It follows from $p_{1,t}=u\circ p_1$ that $\frac{\dw}{\dw t}(\nabla p_1)^T=(\nabla p_1)^T(\nabla u\circ p_1)^T$ and that $\frac{\dw}{\dw t}(m\circ p_1)=(m_t+\nabla m\cdot u)\circ p_1$. Since the determinant is an alternating multilinear form, we see that
$$\frac{\dw}{\dw t}|\nabla p_1|=|\nabla p_1|(\nabla\cdot u)\circ p_1.$$
Consequently we find that
$$\frac{\dw}{\dw t}\{(\rho\circ p_1)|\nabla p_1|\}=\{[\rho_t+\nabla\rho\cdot u+\rho(\nabla\cdot u)]\circ p_1\}|\nabla p_1|=0$$
and, since $p_{2,t}=\rho\circ p_1$,  that
\begin{align}
&\quad\frac{\dw}{\dw t}\{(\nabla p_1)^T(m\circ p_1)|\nabla p_1|+(\nabla p_2)^T(\rho\circ p_1)|\nabla p_1|\}\nonumber\\
&=(\nabla p_1)^T\{[(\nabla u)^Tm+m_t+\nabla m\cdot u+m(\nabla\cdot u)+(\nabla\rho)^T\rho]\circ p_1\}|\nabla p_1|\nonumber\\
&=0.\nonumber
\end{align}
This achieves the proof of our theorem.
\endproof
\section{The case $(\alpha,\beta,\gamma)=(1,0,0)$ and $n=2$}
Let $(\alpha,\beta,\gamma)=(1,0,0)$ and $n=2$, i.e., we consider the $2$D $2$-component $\mu$-Camassa-Holm equation. In this section, it is our aim to compute the sectional curvature of the group $G^\infty$ associated with the $(1,0,0)$-Camassa-Holm equation and to show that this equation is the only member of the so-called multi-variable $b$-equation
\beq\label{beq}m_t=-u\cdot\nabla m-(\nabla u)^T\cdot m-(b-1)m(\nabla\cdot u),\quad m=(\mu-\Delta)u,\eeq
for which the associated geodesic flow is related to a crucial Riemannian structure, as explained in Sect.~\ref{sec_geo}. For any $b>2$, \autoref{beq} belongs to the family of so-called non-metric Euler equations. A similar observation has been made in \cite{K11} for the $(0,1,0)$-Camassa-Holm equation.
Let us now prove the analog of Theorem~3 in \cite{K11}.
\thm Let $b\geq 2$ be an integer and $\mathbb L=\text{\rm diag}(\mu-\Delta,\mu-\Delta)$. Suppose that there is a regular
inertia operator $\A=\text{\rm diag}(A,A)$, $A\in\mathcal L_{\text{\rm is}}^{\text{\rm sym}}(C^{\infty}(\S^2)),$ such that the 2D-$\mu$-$b$-equation
$$m_t=-u\cdot\nabla m-(\nabla u)^Tm-(b-1)m(\nabla\cdot u),\quad m=\mathbb L u,$$
is the Euler equation on $\Diff^{\infty}(\S^2)$ with respect to the right-invariant metric $\rho_\A$ induced by $\A$. Then $b=2$ and $\A=\mathbb L$.\endthm
\endthm\rm
\proof We assume that, for a given $b\geq 2$ and $A\in\mathcal L_{\text{is}}^{\text{sym}}(C^{\infty}(\S^2))$, the
2D-$\mu$-$b$-equation is the Euler equation on the torus
diffeomorphism group with respect to $\rho_\A$. As elucidated in \cite{K11}, this implies that
\bea\label{gl1}\A^{-1}\left\{u\cdot\nabla(\A u)+(\nabla u)^T\A u+\A u(\nabla\cdot u)\right\}&=&\nonumber\\
&&\hspace{-3cm}\mathbb L^{-1}\left\{u\cdot\nabla (\mathbb Lu)+(\nabla u)^T(\mathbb Lu)+(b-1)(\mathbb Lu)(\nabla\cdot u)\right\},\eea
that $\text{span}\{\mathbf 1\}$, $\mathbf 1=(1,1)$, is an invariant subspace of $\A$, more precisely $\A\mathbf 1=\mathbf 1$, and that
\beq\label{gl2}\A^{-1}\left[(\nabla u)^T+\nabla(\A u)+(\nabla\cdot u)\right]\mathbf 1=\mathbb L^{-1}\left[(\nabla u)^T+\nabla(\mathbb L u)+(b-1)(\nabla\cdot u)\right]\mathbf 1.\eeq
Let $n=(n_1,n_2)\in(2\pi\Z)^2\backslash\{(0,0)\}$, and write $z=(x,y)$ for the variable on $\S^2$. We will consider the functions $u_n=e^{\i nz}\mathbf 1$ in the following for which we have $\mathbb Lu_n=n^2u_n$ and $\mathbb L^{-1}u_n=n^{-2}u_n$, $n^2=n_1^2+n_2^2$. Set $v_n=\A u_n$. An explicit calculation of the left hand and right hand side of \eqref{gl2} shows that we have the identity
\beq\label{gl3}\nabla v_n\cdot\mathbf 1-\i\alpha_{n}v_n=-\i\beta_n u_n\eeq
where
\begin{align}
\alpha_n=&\text{ diag}\left((b+1)\frac{n_1}{n^2}+(b-1)\frac{n_2}{n^2}+n_1+n_2,(b+1)\frac{n_2}{n^2}+(b-1)\frac{n_1}{n^2}+n_1+n_2\right),\nonumber\\
\beta_n=&\text{ diag}\left(3n_1+n_2,3n_2+n_1\right).\nonumber\end{align}
Assume that $n_1\neq n_2$ first. Since $v_n=Ae^{\i nz}\mathbf 1$ we get that $Ae^{\i nz}=n^2e^{\i nz}$. For $n_1=n_2$ we see that the function $v_n=\frac{2}{b}n^2e^{\i nz}\mathbf 1$ constitutes a solution to \eqref{gl3}. We insert $u_n$ for $n_1=n_2\neq 0$ into \autoref{gl1} to get that $b=2$. Since $\{e^{\i nz};\,n\in(2\pi\Z)^2\}$ is a basis for $C^{\infty}(\S^2)$, it follows that $\A=\mathbb L$. This completes the proof of the theorem.
\endproof
We also note the following immediate consequence of the above theorem.
\cor The geodesic flow for the two-dimensional and two-component $\mu$-Degasperis-Procesi equation
$$m_t=-u\cdot\nabla m-(\nabla u)^T\cdot m-2m(\nabla\cdot u),\quad m=(\mu-\Delta)u,$$
on the diffeomorphism group of the torus is not related to any right-invariant metric on $\Diff^{\infty}(\S^2)$ with inertia operator $\text{\rm diag}(A,A)$, $A\in\L_{\text{\rm is}}^{\text{\rm sym}}(C^{\infty}(\S^2))$.
\endcor\rm

The existence of a smooth connection
$\bar\nabla$ on a Banach manifold $M$ immediately implies the existence of
a smooth curvature tensor $R$ defined by
$$R(X,Y)Z=\bar\nabla_X\bar\nabla_YZ-\bar\nabla_Y\bar\nabla_XZ-\bar\nabla_{[X,Y]}Z,$$
where $X,Y,Z$ are vector fields on $M$. In
the case of the $(1,0,0)$-Camassa-Holm equation, since there exists a metric $\ska{\cdot}{\cdot}$, we can also define an (unnormalized) sectional curvature
$S$ by
$$S(X,Y) := \langle R(X,Y)Y, X \rangle.$$
In this section, we will derive a convenient formula for $S$ and
use it to determine large subspaces of positive curvature for the
2D 2-component $\mu$-Camassa-Holm equation. Similar results on the $(0,1,0)$-Camassa-Holm equation can be found in \cite{K11}. Let us now state the analog of Theorem~4 of \cite{K11}.
\begin{prop}\label{prop_curv1} Let $R$ be the curvature tensor on
$G^\infty$ associated with the $(1,0,0)$-Camassa-Holm equation. Then $S(u,v) :=\ska{R(u,v)v}{u}$
is given at the identity by
\begin{equation}\label{Suvexpression}
  S(u,v)=\ska{\Gamma(u,v)}{\Gamma(u,v)}-\ska{\Gamma(u,u)}{\Gamma(v,v)}+R(u,v),
\end{equation}
where
\bea R(u,v)\!\!\!&=&\!\!\!\ska{\nabla u\cdot u}{\nabla v\cdot v}-\ska{\nabla u\cdot v}{\nabla u\cdot v}+\ska{\nabla v\cdot u}{\nabla u\cdot v}-\ska{\nabla v\cdot u}{\nabla v\cdot u}\nonumber\\
&&+\ska{[\nabla(\nabla u\cdot u)]\cdot v}{v}-\ska{[\nabla(\nabla u\cdot v)]\cdot v}{u}+\ska{[\nabla(\nabla v\cdot u)]\cdot v}{u}\nonumber\\
&&-\ska{[\nabla(\nabla v\cdot u)]\cdot u}{v}-\ska{\nabla v(\nabla u\cdot u)}{v}-\ska{\nabla u(\nabla v\cdot v)}{u}\nonumber\\
&&+\ska{\nabla v(\nabla v\cdot u)}{u}+\ska{\nabla u(\nabla v\cdot u)}{v}.\nonumber\eea
\end{prop}
\proof
One can use the arguments in the proof of Theorem~4 of \cite{K11}, since the calculations there do not depend on the particular choice of the inertia operator $\A$.
\endproof
Curvature computations have a long tradition in the geometric theory of partial differential equations, see, e.g., \cite{F88,McK82}, and we will now contribute some further examples concerning the $(1,0,0)$-Camassa-Holm equation. For our purposes, two-dimensional subspaces on which $S>0$ are of particular interest since the positivity of $S$ is related to stability properties of the geodesic flow, cf.~\cite{A89}. Large subspaces of positive curvature for variations of the Camassa-Holm equation have been found in \cite{M98,M02,K11,K11b} by evaluating the sectional curvature on trigonometric functions. Here we extend this discussion for the $(1,0,0)$-Camassa-Holm equation.
\exmp Let $\{e_1,e_2\}$ be the canonical basis of $\R^2$, let $k_1,k_2\in 2\pi\N$ and let
$$
v=
\left(
    \begin{array}{c}
      \sin (k_1x)\sin (k_2y) \\
      \sin (k_1x)\sin (k_2y)\\
    \end{array}
\right).
$$
In the following, we will make use of the identity
$$(\mu-\Delta)^{-1}\sin(\alpha x)\cos(\beta y)=\frac{\sin(\alpha x)\cos(\beta y)}{\alpha^2+\beta^2}$$
and the trigonometric formulas
$$\int_\S\cos^2(\alpha x)\dx=\int_\S\sin^2(\alpha x)\dx=\frac{1}{2},\quad\int_\S\cos(\alpha x)\sin(\beta x)\dx=0,\quad\forall\alpha,\beta\in 2\pi\N.$$
We claim that $S(e_i,v)>0$ for $i=1,2$. First we observe that for general $w\in C^{\infty}(\S^2)$, we have that $R(e_i,w)=0$, $i=1,2$; this follows from a short computation using integration by parts. Hence
$$S(e_i,v)=\ska{\Gamma(e_i,v)}{\Gamma(e_i,v)}.$$
We leave it to the reader to perform explicitly the easy calculations leading to the formulas
\bea
S(e_1,v)&=&\frac{1}{8}\frac{2k_1^2+k_2^2}{k_1^2+k_2^2},\nonumber\\
S(e_2,v)&=&\frac{1}{8}\frac{2k_2^2+k_1^2}{k_1^2+k_2^2}.\nonumber
\eea
It follows that, for any $k_1,k_2\in 2\pi\N$, $S>0$ on the spaces $\text{span}\{e_i,v\}$, $i=1,2$.
\endexmp

\end{document}